\documentclass[preprint,number,reqno]{amsart}
\usepackage{doi}
\usepackage{hyperref}
\usepackage[foot]{amsaddr}

\usepackage{amssymb}
\usepackage{amsfonts}
\usepackage{amsmath,amsthm}
\DeclareMathOperator{\diam}{diam}
\DeclareMathOperator{\area}{area}
\newtheorem{theorem}{Theorem}

\newtheorem{question}[theorem]{Question}
\numberwithin{theorem}{section}

\newtheorem{lemma}[theorem]{Lemma}

\title[Remarks on the Crouzeix--Palencia proof]{Remarks on the Crouzeix--Palencia proof that 
the numerical range is a $(1+\sqrt2)$-spectral set}
\thanks{TR supported by grants from NSERC and the Canada Research Chairs Program. \\\indent
FLS supported by grants of the DAAD and the Deutsche Forschungsgemeinschaft (DFG)}

\author{Thomas Ransford}
\address{ D\'epartement de math\'ematiques et de statistique, 
  Universit\'e Laval, Qu\'ebec (QC), G1V 0A6, Canada, }
 \email{ransford@mat.ulaval.ca}
\author{Felix L.~Schwenninger}
\address{Department of  Mathematics, Universit\"at Hamburg, 
  Bundesstr.~55, 20146 Hamburg, Germany}
  \email{felix.schwenninger@uni-hamburg.de}

\begin{document}

%\maketitle
\date{Dec 17, 2017}

\begin{abstract}
Crouzeix and Palencia recently showed that 
the numerical range of a Hilbert-space operator
is a $(1+\sqrt2)$-spectral set for the operator. 
One of the principal  ingredients of their proof can
be formulated as an abstract functional-analysis lemma. 
We give a new short proof of the lemma and show that,
in the context of this lemma, the constant $(1+\sqrt2)$ is sharp.
 \end{abstract}



\maketitle

\section{Introduction}\label{S:intro}

Let $H$ be a complex Hilbert space 
and let $T$ be a bounded linear operator on $H$.
The \emph{numerical range} $W(T)$ of $T$ is defined by
\[
W(T):=\{\langle Tx,x\rangle: x\in H, \|x\|=1\}.
\]
It is a bounded convex set, and is compact if $\dim H<\infty$.

In the recent paper \cite{CP17},  Crouzeix and Palencia,
improving earlier results of Delyon--Delyon \cite{DD99} and Crouzeix \cite{Cr07},
showed that $\overline{W(T)}$ is always a $(1+\sqrt2)$-spectral set for $T$.
This means that, 
for every function $f$ holomorphic on an open set containing $\overline{W(T)}$, 
the operator norm of $f(T)$ satisfies
\begin{equation}\label{E:1+sqrt2}
\|f(T)\|\le (1+\sqrt2)\sup_{z\in W(T)}|f(z)|.
\end{equation}
Crouzeix \cite{Cr04} has conjectured that $(1+\sqrt2)$ may be replaced by $2$.
Simple examples show that the constant $2$ is best possible. 

The point of departure in  all three papers \cite{Cr07,CP17,DD99} is the same.
Let us fix a smoothly bounded open convex set $\Omega$ containing $\overline{W(T)}$.
Denote by $A(\Omega)$ the algebra of continuous functions $f$ on $\overline{\Omega}$
that are holomorphic on $\Omega$, and write $\|f\|_\Omega:=\sup_\Omega|f|$.
As remarked in \cite{DD99},
the containment $\overline{W(T)}\subset\Omega$ 
is reflected in the fact that the operator-valued measure on $\partial\Omega$
\[
\frac{1}{2\pi i}(\zeta I-T)^{-1}\,d\zeta,
\]
used in defining $f(T)$, has positive real part. 
This quickly  leads to the estimate
\begin{equation}\label{E:f+g*}
\|f(T)+(C\overline{f})(T)^*\|\le 2\|f\|_\Omega 
\quad(f\in A(\Omega)),
\end{equation}
where $C\overline{f}$ denotes the Cauchy transform of $\overline{f}$, namely
\[
(C\overline{f})(z):=\frac{1}{2\pi i}\int_{\partial\Omega}\frac{\overline{f(\zeta)}}{\zeta-z}\,d\zeta
\quad(z\in\Omega).
\]
The problem is how to get from \eqref{E:f+g*}, 
which is an estimate on $\|f(T)+(C\overline{f})(T)^*\|$,
to an estimate on $\|f(T)\|$ alone. In \cite{Cr07,CP17,DD99}, 
this is achieved in three different ways.

In \cite{DD99}, the authors prove the invertibility of the map 
$f\mapsto f+\overline{(C\overline{f})}|_{\partial\Omega}$,
considered as a self-map of the space of continuous functions on $\partial\Omega$.
They further obtain  an estimate for the norm of the inverse, which,
together with \eqref{E:f+g*},
leads to the bound
\[
\|f(T)\|\le \Bigl(\Bigl(\frac{2\pi\diam^2(\Omega)}{\area(\Omega)}\Bigr)^3+3\Bigr)\|f\|_\Omega 
\quad(f\in A(\Omega)).
\]

In \cite{Cr07}, Crouzeix estimates $\|(C\overline{f})(T)\|$ directly,
still under the assumption that $\overline{W(T)}\subset\Omega$,
and then uses \eqref{E:f+g*} and the triangle inequality to obtain the  bound
\begin{equation}\label{E:Crouzeix}
\|f(T)\|\le 11.08\|f\|_\Omega \quad(f\in A(\Omega)).
\end{equation}
This bound is universal, in the sense that the constant 11.08 is independent of $\Omega$.
For certain sets $\Omega$, however,  the bound can be improved.
This is discussed in detail in \cite{BCD06}. In particular, if $\Omega$ is a disk,
then $\|f(T)\|\le 2\|f\|_\Omega$, and this is easily seen to be sharp. 
(There are now several proofs of this last inequality, originally
due to Okubo and Ando. 
A particularly short one,  obtained as a simple consequence of \eqref{E:f+g*}, 
can be found in the recent
preprint of Caldwell, Greenbaum and Li \cite{CGL17}.)

Crouzeix's proof of \eqref{E:Crouzeix} is technical and requires a lot of work.
In the Crouzeix--Palencia article \cite{CP17}, 
the passage from \eqref{E:f+g*} to a bound for $\|f(T)\|$ is quite different and much simpler. 
It is effected using an abstract functional-analysis argument, 
which, for convenience, we summarize in the form of a lemma.

\begin{lemma}\label{L:falemma}
Let $T$ be a bounded Hilbert-space operator 
and let $\Omega$ be a bounded open set containing the spectrum of $T$.
Suppose that, for each $f\in A(\Omega)$, 
there exists $g\in A(\Omega)$ such that
\begin{equation}\label{E:conditions}
\|g\|_\Omega\le\|f\|_\Omega
\quad\text{and}\quad
\|f(T)+g(T)^*\|\le 2\|f\|_\Omega.
\end{equation}
Then
\begin{equation}\label{E:falemma}
\|f(T)\|\le(1+\sqrt{2})\|f\|_\Omega
\quad(f\in A(\Omega)).
\end{equation}
\end{lemma}

To prove \eqref{E:1+sqrt2}, 
this lemma is applied with $\Omega$ a smoothly bounded open convex 
set containing $\overline{W(T)}$ and with $g:=C\overline{f}$.
The second inequality in \eqref{E:conditions} 
is then just  \eqref{E:f+g*}.
The first inequality is a fundamental property of the Cauchy transform,
namely that $f\mapsto C\overline{f}$ is a contraction of $A(\Omega)$ into itself 
when $\Omega$ is convex
(see e.g.\ \cite[Lemma~2.1]{CP17}).
Thus \eqref{E:falemma} holds, and the  main result \eqref{E:1+sqrt2} 
then follows upon `shrinking' $\Omega$ down to $W(T)$.

It is quite striking that neither the numerical range nor the Cauchy transform 
appear explicitly Lemma~\ref{L:falemma}.
They enter merely through the inequalities \eqref{E:conditions}.

Our purpose in this note is to give a very short proof of Lemma~\ref{L:falemma}
(even simpler than the argument given in \cite{CP17}),
and to show that, in the context of this lemma, the constant $(1+\sqrt2)$ is sharp.
We conclude the article with a brief discussion 
of how these remarks relate to Crouzeix's conjecture.

\section{Short proof of Lemma~\ref{L:falemma}}\label{S:proof}
Let $K$ denote the norm of the continuous homomorphism 
$f\mapsto f(T):A(\Omega)\to B(H)$.
Our goal is to show that $K\le 1+\sqrt2$.

Let $f\in A(\Omega)$ with $\|f\|_\Omega\le1$. 
By  \eqref{E:conditions}, there exists $g\in A(\Omega)$
such that $\|g\|_\Omega\le 1$ and $\|f(T)+g(T)^*\|\le2$. 
We then have
\[
f(T)f(T)^*f(T)f(T)^*=f(T)\bigl(f(T)+g(T)^*\bigr)^*f(T)f(T)^*-(fgf)(T)f(T)^*,
\]
and, since $fgf\in A(\Omega)$, it therefore follows that
\[
\|f(T)\|^4\le \|f(T)+g(T)^*\|\|f(T)\|^3+\|(fgf)(T)\|\|f(T)\|\le 2K^3+K^2.
\]
Taking the supremum over all  $f$ with $\|f\|_\Omega\le1$, 
we deduce that $K^4\le 2K^3+K^2$, whence $K\le 1+\sqrt2$,
as desired.

\section{Sharpness of Lemma~\ref{L:falemma}}\label{S:sharpness}
We exhibit a pair $(T,\Omega)$, 
satisfying  all the hypotheses  of Lemma~\ref{L:falemma},
such that equality holds in \eqref{E:falemma} 
for a particular (non-zero) choice of $f\in A(\Omega)$.
This shows that, in the context of Lemma~\ref{L:falemma}, 
the constant $(1+\sqrt2)$ is sharp.

Let
\[
T:=\begin{pmatrix}
1 &1\\0&0
\end{pmatrix},
\]
and let $\Omega:=\Omega_0\cup\Omega_1$, 
where $\Omega_j$ is the open disk with center $j$ and radius $1/4$.
Clearly $\Omega$ contains the spectrum of $T$. 
Also, if $f\in A(\Omega)$, then
\begin{equation}\label{E:fidentity}
f(T)=
\begin{pmatrix}
f(1) &f(1)-f(0)\\
0 &f(0)
\end{pmatrix}.
\end{equation}
Indeed, this is clear if $f(z)=z^n$, since $T^n=T$ for all $n\ge1$.
The identity for general $f$ follows 
by the linearity and continuity of the map $f\mapsto f(T)$.
(It could also be deduced directly from the definition of $f(T)$ as a Cauchy integral.)

Given $f\in A(\Omega)$, define $g\in A(\Omega)$ by
\[
g(z):=
\begin{cases}
-\overline{f(0)}, &z\in\Omega_0,\\
-\overline{f(1)}, &z\in\Omega_1.
\end{cases}
\]
Clearly $\|g\|_\Omega\le\|f\|_\Omega$. Also
\begin{align*}
\|f(T)+g(T)^*\|
&=\Bigl\|
\begin{pmatrix}
f(1) &f(1)-f(0)\\
0 &f(0)
\end{pmatrix}
+
\begin{pmatrix}
-f(1) &0\\
-(f(1)-f(0)) &-f(0)
\end{pmatrix}
\Bigr\|\\
&=|f(1)-f(0)|\,\Bigl\|
\begin{pmatrix}
\hfill0 &1\\ -1 &0
\end{pmatrix}\Bigr\|
\le2\|f\|_\Omega.
\end{align*}
Thus \eqref{E:conditions} holds.

Finally, let $h\in A(\Omega)$ be given by
\[
h(z):=
\begin{cases}
-1, &z\in\Omega_0,\\
\hfill1, &z\in\Omega_1.
\end{cases}
\]
Then $\|h\|_\Omega=1$ and
\[
\|h(T)\|
=\Bigl\|
\begin{pmatrix} 
h(1) &h(1)-h(0)\\
0 &h(0)
\end{pmatrix}
\Bigr\|
=\Bigl\|
\begin{pmatrix} 
1 &\hfill2\\
0 &-1
\end{pmatrix}
\Bigr\|=1+\sqrt2.
\]
Thus equality holds in \eqref{E:falemma} when $f=h$.

\section{Crouzeix's conjecture}
As mentioned in the introduction, 
Crouzeix has conjectured that 
$\overline{W(T)}$ is always a $2$-spectral set for $T$. 
What does the example in \S\ref{S:sharpness} tell us about the conjecture?

First of all, the example is \emph{not} a counterexample to the conjecture, 
because in the example $W(T)\not\subset\Omega$.
Indeed, the eigenvalues $0,1$ of $T$ 
belong to different components $\Omega_0,\Omega_1$ of $\Omega$,
whereas the numerical range $W(T)$ is a convex set  containing $0$ and $1$. 
In fact, Crouzeix's conjecture is known to be true 
for $2\times2$ matrices (see \cite{Cr04}).

What the example \emph{does} tell us 
is that it is not possible to improve upon the constant $(1+\sqrt2)$ merely
by adjusting the proof of Lemma~\ref{L:falemma}. 
In this sense, the grouping of terms in the proof presented in \S\ref{S:proof} is already optimal. 

If one is to approach the conjecture along the same lines as in Lemma~\ref{L:falemma},
then more information is needed about the choice of $g$. 
Here is one possibility. 
Recall that, in the application of Lemma~\ref{L:falemma},
$g=C\overline{f}$, the Cauchy transform of $\overline{f}$. 
The map $f\mapsto C\overline{f}$ is both antilinear and unital (i.e. it maps $1$ to $1$). 
In the example in \S\ref{S:sharpness}, 
the map $f\mapsto g$ is also antilinear, 
but it is \emph{not} unital; on the contrary it sends $1\mapsto -1$. 
This suggests the following question.

\begin{question}\label{Q:faq}
Let $T$ be a bounded Hilbert-space operator 
and let $\Omega$ be a bounded open set containing the spectrum of $T$.
Suppose that there exists a unital antilinear  map 
$\alpha:A(\Omega)\to A(\Omega)$ such that,
for all $f\in A(\Omega)$,
\[
\|\alpha(f)\|_\Omega\le\|f\|_\Omega
\quad\text{and}\quad
\|f(T)+(\alpha(f))(T)^*\|\le 2\|f\|_\Omega.
\]
Does it follow that
\[
\|f(T)\|\le 2\|f\|_\Omega \quad(f\in A(\Omega))?
\]
\end{question}

An affirmative answer to this question would prove the Crouzeix conjecture.

\section*{Acknowledgments}
The results in this note were obtained at the Workshop on Crouzeix's Conjecture
held at the American Institute of Mathematics (AIM). The authors express their
gratitude both to AIM  and to the organizers of the workshop.

%\bibliographystyle{siamplain}
%\bibliography{references}

\end{document}